\documentclass[oneside, 11pt, reqno]{amsart}
\input{preamble.sty}
\title{Characterising canonical lifts of  ordinary Calabi--Yau varieties}
\author{Samuel A.\ Moore}
\address{St. John's College, University of Oxford}
\email{samuel.moore@maths.ox.ac.uk}
\usetikzlibrary{decorations.pathmorphing}

\newcommand{\syn}{\mrm{syn}}
\newtheorem{thmx}{Theorem}

\begin{document}

\begin{abstract} We study Brantner--Taelman's canonical lifts of ordinary Calabi--Yau varieties and characterise their (quasi)canonical lifts in terms of rational $p$-adic Hodge theoretic data. We also show that these canonical lifts need not be of CM type and that the associated local Galois representations need not split as sums of twists of unramified representations, in contrast to classical results for abelian varieties and K3 surfaces.
\end{abstract}

\maketitle

\section{Introduction}

Fix a prime $p$ and a perfect field $k$ of characteristic $p$. By a celebrated theorem of Serre--Tate, any ordinary abelian variety $Y$ over $k$ admits a canonical lift to the Witt vectors $W(k)$ \cite{ST-original} \cite{Katz-ST}. More precisely, Serre and Tate constructed a formal torus structure on the deformation functor of $Y$, and the canonical lift is defined to correspond to the identity. An analogous theory for K3 surfaces was developed by Deligne--Illusie and Nygaard \cite[\S 2.1 (B)]{DI-can-coords} \cite[\S 1]{NygaardTC}, and there has been work extending Serre--Tate theory to wider classes of Calabi--Yau varieties, e.g.\ \cite{Ward}, \cite{AZ-ST}.

Most recently, Brantner and Taelman have constructed a canonical equivalence $\mrm{st}_Y$ between the derived deformation functor $\mrm{Def}_Y$ of a Bloch--Kato $2$-ordinary Calabi--Yau variety $Y$ over a perfect field $k$ (satisfying a torsionfree condition), and a derived formal group $\mrm{ST}_Y$ of multiplicative type \cite[Theorem B]{BT}. In particular, any such $Y$ admits a canonical lift to $\mrm{Spf}\,W(k)$, corresponding to the identity section; if $Y$ is projective, it is algebraisable.

In this note, we study Brantner--Taelman's canonical lifts and the more general \emph{quasicanonical} lifts:

\begin{definition} (cf. \cite[p.\ 7]{ST-original}) Let $\mathcal O_K$ be a complete mixed characteristic DVR with residue field $k$. A \emph{quasicanonical lift} of $Y$ is a formal lift $\mathfrak X/\Spf(\mathcal O_K)$ which corresponds under Brantner--Taelman's isomorphism $\mrm{st}_Y \colon \pi_0 \mrm{Def}_Y(\mrm{Spf}(\mathcal O_K)) \xrightarrow{\cong} \pi_0 \mrm{ST}_Y(\Spf(\mathcal O_K))$ to a torsion point of $\pi_0 \mrm{ST}_Y(\mrm{Spf}(\mathcal O_K))$.
\end{definition}
Here we put $F(\mrm{Spf}(\mathcal O_K)) = \varprojlim_n F(\mathcal O_K/\mathfrak m^n)$.
\begin{remark}  Under the torsionfree conditions of \cref{thmA} below, the distinction between quasicanonical and canonical lifts disappears when $p-1$ does not divide the ramification index $e_K$; see \cref{simplify}.
\end{remark}

\subsection{Statement of results}

Our first result concerns characterising quasicanonical lifts in terms of familiar $p$-adic invariants. To state the result, let $K$ be a complete discretely valued field of characteristic $0$ with ring of integers $\mathcal O_K$ and perfect residue field $k$ of characteristic $p > 0$, and put $K_0 = W(k)[1/p] \subset K$. Fix an algebraic closure $\overline{K}$ with completion $C$.

\begin{thmx} \label{thmA} 
Let $\mathfrak X_0$ be a Bloch--Kato $2$-ordinary Calabi--Yau variety over $k$ of dimension $d$, and put $\overline{\mathfrak X}_0 = \mathfrak X_0 \times_k \overline{k}$.  Suppose that
  \begin{enumerate} \item[(i)] The $\Z_p$-modules $\mrm{H}^d_{\et}(\overline{\mathfrak X}_0, \G_m^{\land})$ and $\mrm{H}^d_{\et}(\overline{\mathfrak X}_0, \Z_p)$ are torsionfree, and
  \item[(ii)] The absolute ramification index $e_K$ is less than $p$.
  \end{enumerate}
 Let $\mathfrak X/\mrm{Spf}(\mathcal O_K)$ be a flat formal lift of $\mathfrak X_0$. Then the following are equivalent:
\begin{enumerate} \item $\mathfrak X$ is a quasicanonical lift of $\mathfrak X_0$.
\item \label{2} The Newton slope $1$ subspace $\mrm{H}^d_{\cris}(\mathfrak X_0/K_0)^{(1)} \subseteq \mrm{H}^d_{\cris}(\mathfrak X_0/K_0)$ is contained in $\mrm{Fil}^1 \mrm{H}^d_{\dR}(\mathfrak X_K/K)$, where we identify $\mrm{H}^d_{\cris}(\mathfrak X_0/K_0) \otimes_{K_0} K \simeq \mrm{H}^d_{\dR}(\mathfrak X_K/K)$.
\item \label{3} There exists a subrepresentation $V(-1) \subseteq \mrm{H}^d_{\et}(\mathfrak X_{C}, \Q_p)$ such that $V$ is unramified and $1$ is not a Hodge--Tate weight of the quotient.
\end{enumerate}
\end{thmx}

See $(\star)$ below for a comparison to analogous characterisations of quasicanonical lifts of ordinary abelian varieties and K3 surfaces.

\begin{remark}
  \begin{enumerate}
    \item In \hyperref[2]{(\ref*{2})},  the $(-1)$ denotes a Tate twist. In \hyperref[3]{(\ref*{3})}, we refer to the (Newton) slope $1$ subspace of the isocrystal $\mrm{H}^d_{\cris}(\mathfrak X_0/K_0)$, in the sense of the Dieudonn\'{e}--Manin classification. 
    \item The Bloch--Kato $2$-ordinarity and torsion conditions are implied by the Newton and Hodge slope $0$ and $1$ multiplicities of the isocrystals $\mrm{H}^i_{\cris}(\mathfrak X_0/K_0)$ agreeing, i.e.\ that for $r = 0,1$ and each $i \ge 0$, $\dim_{\Q_p} (\mrm{H}^i_{\cris}(\mathfrak X_0/K_0) \otimes_{K_0} \widehat{K_0^{\mrm{unr}}})^{\phi = p^r} = \dim_{k} \mrm{H}^{i-r}(\mathfrak X_0, \Omega_{\mathfrak X_0/k}^r)$ \cite[Proposition 7.3 (6)]{BlochKato}. 
    \end{enumerate}
  \end{remark}

  \begin{remark}
    The unramified representation $V$ in \hyperref[3]{condition (3)} above is necessarily given by $\mrm{H}^d_{\et}(\mathfrak X_{C}, \Q_p(1))^{I_K}$. The condition can therefore be rephrased more intrinsically as the usual inequality $\dim_{\Q_p} \mrm{H}^d_{\et}(\mathfrak X_C, \Q_p(1))^{I_k} \le \dim_{C} \mrm{H}^{d-1}(\mathfrak X_K, \Omega^1_{\mathfrak X_K/K})$ being an equality. \end{remark}

\cref{addendum} gives an application to derived categories, showing that often we can lift derived equivalences to canonical lifts \cref{der-equiv}. This is closely related to work of Kurama \cite[Theorem 1.1, Question 1.9]{kuramaFM}.
  
Our second main result concerns the general properties of (quasi)canonical lifts of ordinary Calabi--Yau varieties. The quasicanonical lifts of ordinary abelian varieties and K3 surfaces over finite fields are known to be of complex multiplication type (henceforth \emph{CM-type}).  For the purposes of this paper, being of CM-type means that the $\ell$-adic algebraic monodromy groups $G_\ell(V)$ are abelian for each prime $\ell$; see \cref{cm}. Moreover, one can characterise canonical (or quasicanonical) lifts of a K3 surface or abelian variety among lifts to $\mrm{Spec}(\mathcal O_K)$ by the following property:
\begin{enumerate}
\item[($\star$)] The $G_{K}$-representation $\mrm{H}^d_{\et}(X_{\overline{K}},\Z_p)$ (resp.\ $\mrm{H}^d_{\et}(X_{\overline{K}}, \Q_p)$) splits as a direct sum of Tate twists of unramified representations, where $d = \dim X$.
\end{enumerate}
For abelian varieties, this is true more or less by definition (and fully-faithfulness of the Tate module construction), and it is \cite[Theorem 2.6]{NygaardTC} \cite[Lemma 1.5.]{JDY} \cite[Theorem C]{TaelmanOrdinaryK3} for K3 surfaces.

We show that neither of these results generalises to canonical lifts of Calabi--Yau threefolds by studying the good mod $p$ reductions of a rigid Calabi--Yau threefold $V$ over $\Q$ studied by Saito--Yui \cite{Saito-Yui}, which has a smooth projective model $X \to \Spec \Z[1/6]$.
\begin{thmx} \label{thmB}
The Calabi--Yau threefold $V$ over $\Q$ (and its integral model $X$) satisfy the following properties:
\begin{enumerate}
\item For any prime $\ell$, the $\ell$-adic algebraic monodromy group $G_\ell(V)$ is isomorphic to $\mrm{GL}_{2, \Q_\ell}$; in particular $V$ is not of ($\ell$-adic) CM-type.
  \item If $5 \le p < 47351$, the mod $p$ reduction $Y_p = X_{\F_p}$ is an ordinary, rigid, strict Calabi--Yau threefold with torsionfree crystalline cohomology. In particular, $X_p = X_{\Spec \Z_p}$ is the canonical lift of $Y_p$ to $\Z_p$.
\item The two-dimensional $G_{\Q_7}$-representation $\mrm{H}^3_{\et}(V_{\overline{\Q}}, \Q_7)$ is indecomposable.
\end{enumerate}
\end{thmx}

This example also shows that some feasible strengthenings of conditions (2) and (3) in \cref{thmA} fail in general for canonical lifts; see \cref{strength}.

\subsection{Acknowledgments}

It is a pleasure to thank our advisor, Lukas Brantner, for all his support and for suggesting the questions studied here. We also thank Simon Felten for helpful conversations about threefolds.

\subsection{AI disclosure}

Sections \S 1-4 of this paper were developed without AI tools beyond proofreading, with the exception that an LLM made us aware of \cite{yyang26}, leading to removing a torsion hypothesis in \cref{lemm}; see there for details. \cref{noncm} benefitted more substantially from LLMs: whilst the idea of considering the mod $p$ reduction of a rigid modular CY3 $V$) was ours, an LLM proposed studying this particular $V$ and using work of Schoen. Finally, during proofreading, we were reminded of the derived invariance of the Mukai structure by an LLM; together with our results, this led quickly to the main result of \cref{addendum}.

\subsection{Terminology and conventions} We mostly follow the notation of \cite{BT}. For the purposes of this paper, a Calabi--Yau variety over a field $F$ is a geometrically irreducible smooth and proper $F$-scheme with trivial canonical bundle; a \emph{strict} Calabi--Yau variety $Y$ additionally satisfies $\mrm{H}^i(Y, \mathcal O_Y) = 0$ for $0 < i < \dim(Y)$. For \S \S 2-4, we fix $k, K_0, \mathcal O_K, \dots$ as in the introduction, and write $\breve{K} = \widehat{K^{\mathrm{nr}}}$.

\section{On the Serre--Tate domain}

\begin{notation} For \S \S 2--3, we let $Y$ denote a smooth and proper $k$-scheme of pure dimension $d$ and structure morphism $f \colon Y \to \Spec k$. We write $S = \Spec k$.
\end{notation}

\subsection{Recollections}

\begin{definition} For an integer $n \ge 1$, $Y$ is \emph{Bloch--Kato $n$-ordinary} (or simply \emph{BK} $n$-ordinary) if $\mathrm{H}^*(Y, B^r_{Y/k}) = 0$ for all $1 \le r \le n$, where $B^r_{Y/k} = \mrm{im}(d \colon \Omega^{r-1}_{Y/k} \to \Omega^r_{Y/k})$. We call $Y$ \emph{Bloch--Kato ordinary} if it is Bloch--Kato $n$-ordinary for all $n$.
\end{definition}

\begin{remark} \label{bkrecoll} Equivalently, $Y$ is $1$-ordinary if and only if the map $\F_p \to \mathcal O_Y$ of \'{e}tale sheaves induces an equivalence $\mathcal O_{S} \otimes_{\F_p} \mrm{R}f_* \F_p \xrightarrow{\simeq} \mrm{R}f_* \mathcal O_Y$, and, given that it is $1$-ordinary, BK $2$-ordinarity is equivalent to the map $\mrm{dlog} \colon \mathbf G_{m,Y}^{\land} \to \Omega^1_{Y/k}$ inducing an equivalence $\mathcal O_S \otimes_{\Z_p} \mrm{R}f_* \mathbf G_{m,Y}^{\land} \xrightarrow{\simeq} \mrm{R}f_* \Omega^1_{Y/k}$ \cite[\S 7]{BlochKato}, \cite[Proposition 8.11]{BT}.
\end{remark}

We recall that Brantner--Taelman's Serre--Tate domain is defined on animated artinian $W(k)$-algebras with residue field $k$ as 
\[ \mrm{ST}_Y(-) = \mrm{Map}_{\mrm{Mod}(S_{\et}, \Z_p)}(\mrm{R}f_* \G_m^{\land},\, \mrm{R}^d f_* \Z_p \otimes_{\Z_p} \mathbf G_m^{\land}(\mathcal W \otimes_W (-))[1-d]).\]
Here we write $\mathcal W$ for the sheaf of Witt vectors on the small \'{e}tale site $S_{\et}$, and $\mathbf G_m^{\land}(-)$ denotes the (sheafified) $p$-completed units.

\begin{situation} \label{sit} Let $Y$ be a BK $2$-ordinary Calabi--Yau variety over $k$. We say it satisfies condition $(\dagger$) if the $\Z_p$-modules $\mathrm{H}^d_{\et}(Y_{\overline{k}}, \Z_p)$ and $\mrm{H}^d_{\et}(Y_{\overline{k}}, \G_m^{\land})$ are torsionfree. 
\end{situation}

\subsection{A simplified description}

We wish to give a slightly simpler description of $\pi_0 \mrm{ST}_Y$ under condition ($\dagger$) above. The key cohomological input is the following:
\begin{lemma} \label{connective} Let $Y$ be a BK $2$-ordinary Calabi--Yau variety satisfying condition $(\dagger)$. For any artinian local $W(k)$-algebra $A \in \mrm{CRing}_{W}^{\mrm{art}}$, the mapping spectrum \[G(A) =  \mrm{map}_{\Mod(S_{\et}, \Z_p)}(\mrm{R}^d f_* \G_m^{\land}, \mrm{R}^d f_* \Z_p \otimes_{\Z_p} \G_m^{\land}(\mathcal W \otimes_W A))\] is discrete.
\end{lemma}

\begin{proof}
We induct on the length of $A$. If $A = k$, then $G(A)$ vanishes; otherwise, we can write $A$ as a square-zero extension of some $B \in \mrm{CRing}^{\mrm{art}}_{W}$ by $k$; writing $A = B \times_{\mrm{sqz}_k(k[1])} k$ for some map $k \to \mrm{sqz}_k(k[1])$, we obtain a fibre sequence as follows (cf.\ the proof of \cite[Lemma 8.15]{BT}):
  \[ \mathcal O_{S} \to \G_m^{\land}(\mathcal W \otimes_W A) \to \G_m^{\land}(\mathcal W \otimes_W B).\]
  Applying $\mrm{map}_{\mrm{Mod}(S_{\et}, \Z_p)}(\mrm{R}^d f_* \G_m^{\land}, \mrm{R}^d f_* \Z_p \otimes_{\Z_p}(-))$, we obtain a fibre sequence \[K = \mrm{map}_{\mrm{Mod}(S_{\et}, \Z_p)}(\mrm{R}^d f_* \G_m^{\land}, \mrm{R}^{d} f_* \Z_p \otimes_{\Z_p} \mathcal O_S) \to G(A) \to G(B).\] Inducting on length, it suffices to show that $K$ is discrete. By Bloch--Kato $1$-ordinarity and as $\mrm{H}^d_{\et}(Y_{\overline{k}}, \Z_p)$ is torsionfree, $\mrm{R}^d f_* \Z_p \otimes_{\Z_p} \mathcal O_S \simeq \mrm{R}^d f_* \F_p \otimes_{\F_p} \mathcal O_S \simeq \mrm{R}^d f_* \mathcal O_Y$; similarly $\mrm{R}^d f_* \mathbf G_m^{\land} \otimes_{\Z_p} \mathcal O_S \simeq \mrm{R}^d f_* \Omega^1_{Y/k}$. We can now compute
  \begin{align*}
    K &\simeq \mrm{map}_{\mrm{Mod}(S_{\et}, \Z_p)}(\mrm{R}^d f_* \G_m^{\land}, \mrm{R}^d f_* \mathcal O_Y) \\
    & \simeq \mrm{map}_{\mrm{Mod}(S_{\et}, \mathcal O_S)}(\mrm{R}^{d} f_* \G_m^{\land} \otimes_{\Z_p} \mathcal O_S, \mrm{R}^d f_* \mathcal O_Y) \\
      &\simeq \mrm{map}_{\mrm{QCoh}(S)}(\mrm{R}^d f_* \Omega^1_{Y/k}, \mrm{R}^d f_* \mathcal O_Y)
 \\ &\simeq \mrm{map}_{\Mod_k}(\mrm{H}^d(Y, \Omega^1_{Y/k}), \mrm{H}^d(Y, \mathcal O_Y)), \end{align*} which is discrete.
\end{proof}

\begin{proposition} \label{simplify} Let $Y$ be a BK $2$-ordinary Calabi--Yau variety satisfying $(\dagger)$, and let $A \in \mrm{CRing}^{\mrm{art}}_{W}$.

  \begin{enumerate} \item \label{simp1} Applying $\pi_{1-d}$ induces an isomorphism
    \begin{align*} \pi_0 \mrm{ST}_Y(A) &\xrightarrow{\cong} \mrm{Hom}_{\mrm{Mod}(S_{\et}, \Z_p)}(\mrm{R}^{d-1}f_* \G_m^{\land}, \mrm{R}^d f_* \Z_p \otimes_{\Z_p} \G_m^{\land}(\mathcal W \otimes_W A)) \\ & \cong \mrm{Hom}_{\mrm{Gal}(\overline{k}/k)}(\mrm{H}^{d-1}_{\et}(Y_{\overline{k}}, \G_m^{\land}), \mrm{H}^d_{\et}(Y_{\overline{k}}, \Z_p) \otimes_{\Z_p} (1 + \mathfrak m_{{A}^{\mrm{sh}}})),
    \end{align*}
    where $A^{\mrm{sh}} = A \otimes_{W(k)} W(\overline{k})$. Moreover, given a lift $X/A$, the corresponding map fits in a canonical exact sequence
   \[\label{eq} \mrm{H}^{d-1}_{\et}(\mathfrak X_{A^{\mrm{sh}}}, \G_m^{\land}) \to \mrm{H}^{d-1}_{\et}(Y_{\overline{k}}, \G_m^{\land}) \to \mrm{H}^d_{\et}(Y_{\overline{k}}, \Z_p) \otimes_{\Z_p} (1 + \mathfrak m_{A^{\mrm{sh}}}).
    \]
  \item Suppose $k'/k$ is an extension of perfect fields and put $\mathcal O_L = \mathcal O_K \otimes_{W(k)} W(k')$. Then for any lift $\mathfrak X/\mathcal O_K$ of $Y$, $X$ is a canonical (or quasicanonical) lift of $Y$ to $\mathcal O_K$ if and only if $X' = X \otimes_{\mathcal O_K} \mathcal O_L$ is a (quasi)canonical lift of $Y' = Y \otimes_k k'$ to $\mathcal O_L$.
  \item A formal lift $\mathfrak X/\mrm{Spf}(\mathcal O_K)$ is a canonical (or quasicanonical) lift if and only if the map $\mrm{H}^{d-1}(\mathfrak X_{\mrm{Spf}(\breve{\mathcal O}_K)}, \G_m^{\land}) \to \mrm{H}^{d-1}_{\et}(Y_{\overline{k}}, \G_m^{\land})$ is surjective (resp.\ rationally surjective).
    \item If $p-1 \nmid e_K$, there is only one quasicanonical lift of $Y$ to $\mathcal O_K$, namely the canonical lift.
  \end{enumerate}
\end{proposition}

\begin{proof} For (1), we have in general that \[\pi_0 \mrm{ST}(A) \simeq \Hom_{\Mod(S_{\et}, \Z_p)}(\tau_{\le 1-d} \mrm{R}f_* \G_{m,Y}^{\land}, \mrm{R}^df_* \Z_p \otimes_{\Z_p} \mathbf G_m^{\land}(\mathcal W \otimes_W A)[1-d]).\] Since $\mrm{R}f_* \mathbf G_{m,Y}^{\land}$ is concentrated in homological degrees $[-d,0]$ (e.g.\ by \cref{bkrecoll}), the first isomorphism in part (1) now follows from \cref{connective}; the second is a standard identification.

The exact sequence in part (1) obtains by unwinding how Brantner--Taelman's period map is defined in \cite[\S 8]{BT}. The period map \[\mrm{Def}_Y(A) \xrightarrow{\cong} \mrm{ST}_Y(A) = \mrm{Map}_{\Mod(S_{\et}, \Z_p)}( \mrm{R}f_* \G^{\land}_{m,Y}, \mrm{R}^d f_* \Z_p \otimes_{\Z_p} \mathbf G_m^{\land}(\mathcal W \otimes_W A)[1-d])\] is given on $\pi_0$ by sending a lift $X/A$ to the map $c$ in a cofibre sequence
\begin{equation} \label{eq} \mrm{R}f_* \mathbf G_m^{\land}(\mathcal O_X) \to \mrm{R}f_* \mathbf G_{m}^{\land}(\mathcal O_Y) \xrightarrow{c} \mrm{R}f_* \mathbf G_m^{\land}(\mathcal W \otimes_W A)[1]
\end{equation}
and composing with the truncation map to $\mrm{R}^d f_* \Z_p \otimes_{\Z_p} \mathbf G_m^{\land}(\mathcal W \otimes_W A)[1-d]$. Taking $\pi_{1-d}$, we obtain an exact sequence as in (1).

For part (2), write $f' \colon Y' \to \Spec k'$ for the structure morphism. Any lift $Z$ of $Y$ to $\mathcal O_K$ pulls back along $\Spf(\mathcal O_L) \to \Spf(\mathcal O_K)$ to a lift of $Y'$ to $\mathcal O_L$; this furnishes a map $\Def_{Y}(\mrm{Spf}\,\mathcal O_K) \to \mrm{Def}_{Y'}(\mrm{Spf}\,\mathcal O_L)$ which fits into a commutative diagram as follows:
\[\begin{tikzcd}[column sep = 0.8 em]
    {\pi_0 \mrm{Def}_Y(\mrm{Spf}\, \mathcal O_K)} \arrow[r, "\cong"] \arrow[d] & {\pi_0 \mrm{ST}_Y(\mrm{Spf}\, \mathcal O_K)} \arrow[r, "\cong"] \arrow[d] & {\mrm{Hom}_{G_k}(\mrm{H}^{d-1}_{\et}(Y_{\overline{k}}, \G_m^{\land}), \mrm{H}^d_{\et}(Y_{\overline{k}}, \Z_p) \otimes_{\Z_p} (1 + \mathfrak m_{\mathcal O_K^{\mrm{sh}}}))} \arrow[d] \\
 {\pi_0 \mrm{Def}_{Y'}(\mrm{Spf}\, \mathcal O_L))} \arrow[r, "\cong"] & {\pi_0 \mrm{ST}_{Y'}(\Spf \mathcal O_L)} \arrow[r, "\cong"] & {\mrm{Hom}_{G_{k'}}(\mrm{H}^{d-1}_{\et}(Y_{\overline{k'}}, \G_m^{\land}), \mrm{H}^d_{\et}(Y_{\overline{k'}}, \Z_p) \otimes_{\Z_p} (1 + \mathfrak m_{\mathcal O_L^{\mrm{sh}}}))}
  \end{tikzcd}
 \]
The arrow on the right is injective: the functor from $G_k$- to $G_{k'}$-representations is faithful and the map $\mathcal O_K^{\mrm{sh}} \to \mathcal O_L^{\mrm{sh}}$ is injective. As the canonical (resp.\ a quasicanonical) lift of $Y$ corresponds to the identity element (resp.\ a torsion element) of each hom, part (2) obtains.

Part (3) is essentially formal from (1). Given a formal lift, it is a canonical (resp.\ quasicanonical) lift if and only if the map $\mrm{H}^{d-1}_{\et}(Y_{\overline{k}}, \G_m^{\land}) \to \mrm{H}^{d}_{\et}(Y_{\overline{k}}, \Z_p) \otimes_{\Z_p} (1 + \mathfrak m_{\breve{\mathcal O}_K})$ vanishes (resp.\ has torsion image, or equivalently vanishes rationally). This is equivalent, by passage to the limit of the fibre sequences in \hyperref[eq]{equation (\ref*{eq})} and $p$-adic continuity (cf.\ \cite[Remark 7.2.3]{BL}), to the map $\mrm{H}^{d-1}_{\et}(\mathfrak X_{\breve{\mathcal O}_K}, \G_m^{\land}) \to \mrm{H}^{d-1}_{\et}(Y_{\overline{k}}, \G_m^{\land})$ being (rationally) surjective.

Finally, for part (4), we observe that (under the torsion hypotheses), the existence of torsion in $\pi_0 \mrm{ST}_Y(\mrm{Spf}(\mathcal O_K))$ implies that of (necessarily $p$-)torsion in $1 + \mathfrak m_{\breve{\mathcal O_K}} \subset \widehat{\mathcal O_K}^\times$, i.e a $p$th root of unity in $\breve{K}$. Since $e_{\breve{K}_0(\zeta_p)/\breve{K_0}} = p-1$ and $e_{K/K_0} = e_{\breve{K}/\breve{K}_0}$, we obtain the claimed result. \end{proof}

\section{Characterisation and the Beilinson fibre square}

To provide a simple characterisation of quasicanonical lifts of Calabi--Yau varieties, we will use the following form of a celebrated square of Antieau--Mathew--Morrow--Nikolaus appearing in \cite{AMMN} (whose various antecedents include \cite{Goodwillie}, \cite{DGM}, \cite{Beil14}, \cite{BEK}). Below we write $\mrm{qSyn}_{\Z_p}$ for the category of $p$-torsionfree quasisyntomic rings,  $\mrm{dR}$ for $p$-completed derived de Rham cohomology, and use $(-)_{\Q_p}$ to denote (derived) $p$-completion before rationalisation.

\begin{theorem}\label{var0} \emph{(The Beilinson fibre square on graded pieces \cite[Theorem 6.17]{AMMN})}

Let $R \in \mrm{qSyn}_{\Z_p}$ be a $p$-torsionfree quasisyntomic ring and $i \ge 0$ be an integer. There exists a functorial fibre square in $\Mod_{\Q_p}$ as follows:
\[\begin{tikzcd}
{\RG_{\mrm{syn}}(R, \Q_p(i))} \arrow[r] \arrow[d]  & {\RG_{\mrm{syn}}(R/p, \Q_p(i))} \arrow[d] \\ {(\Fil^i \dR_{R/\Z_p})_{\Q_p}} \arrow[r] & {(\dR_{R/\Z_p})_{\Q_p}}
\end{tikzcd}
\]
\end{theorem}

We will make use of the following variant and identification; this square also appears in the proof of \cite[Theorem 7.13]{AMMN}.

\begin{proposition}\label{var}
Let $\mathfrak X$ be a $p$-torsionfree quasisyntomic formal $\mathcal O_K$-scheme and $i \ge 1$. Using the map $\dR_{\mathfrak X/\Z_p} \to \dR_{\mathfrak X/\mathcal O_K}$, we obtain a fibre square
 \[\begin{tikzcd}
    {\RG_{\mrm{syn}}(\mathfrak X, \Q_p(i))} \arrow[r] \arrow[d] & {\RG_{\mrm{syn}}(\mathfrak X_0, \Q_p(i))} \arrow[d, "\psi"] \\
    {\RG(\mathfrak X, \Fil^{i}\dR_{\mathfrak X/\mathcal O_K})_{\Q_p}} \arrow[r] & {\RG(\mathfrak X, \dR_{\mathfrak X/\mathcal O_K})_{\Q_p}.}
  \end{tikzcd}
\]
If $e_K < p$, then the composite $\RG_{\syn}(\mathfrak X_0, \Q_p(i)) \xrightarrow{\psi} \RG_{\dR}(\mathfrak X/\mathcal O_K)_{\Q_p} \simeq \RG_{\cris}(\mathfrak X_0/W) \otimes_W K$ is naturally homotopic to the composite
\[\RG_{\syn}(\mathfrak X_0, \Q_p(i)) \simeq (\RG_{\cris}(\mathfrak X_0/W)[1/p])^{\phi = p^i}  \to \RG_{\cris}(\mathfrak X_0/W)[1/p] \otimes_{K_0} K.\]
\end{proposition}

\begin{proof}
 By Zariski descent, we may assume that $\mathfrak X = \mrm{Spf}(R)$ and $\mathfrak X_0 = \Spec R_0$ are affine. For any $\mathcal O_K$-algebra $R$, the following square is cartesian,  \[\begin{tikzcd}
      {(\Fil^i \dR_{R/\Z_p})_{\Q_p}}\arrow[r] \arrow[d] & {(\dR_{R/\Z_p})_{\Q_p}} \arrow[d] \\
      {(\Fil^i \dR_{R/\mathcal O_K})_{\Q_p}} \arrow[r] & {(\dR_{R/\mathcal O_K})_{\Q_p}.} 
      \end{tikzcd}
    \]
To prove this, it suffices to show that $(\mrm{L}_{R/\Z_p})_{\Q_p} \to (\mrm{L}_{R/\mathcal O_K})_{\Q_p}$ is an equivalence. By considering the transitivity sequence for $\Z_p \to \mathcal O_K \to R$, it is enough to show that $(\mrm{L}_{\mathcal O_K/\Z_p})_{\Q_p} \simeq 0$. Now $\mrm{L}_{W(k)/\Z_p}^{\land} \simeq 0$, since $\mrm{L}_{W(k)/\Z_p} \otimes_{\Z_p} \F_p \simeq \mrm{L}_{k/\F_p} = 0$, so $\mrm{L}_{\mathcal O_K/\Z_p}^{\land} \simeq \mrm{L}_{\mathcal O_K/W(k)}^{\land}$ by the transitivity sequence. But $\mrm{L}_{\mathcal O_K/W(k)}$ is already $p$-complete, and $\mrm{L}_{\mathcal O_K/W(k)}[1/p] \simeq \mrm{L}_{K/K_0} = 0$. 

Pasting with the square of \cref{var0}, we now obtain the square in \cref{var} except that $R_0$ is replaced by $R/p$; the map $\Q_p(i)(R/p) \to \Q_p(i)(R_0)$ is an equivalence, however, since the Frobenius acts as multiplication by $p^i$ and a power of the Frobenius on $R/p$ factors through $R_0$.

Now there is a functorial commutative diagram
\[\begin{tikzcd}
    {\dR_{R/W}} \arrow[rr] \arrow[d, "\simeq"] & & {\dR_{R/\mathcal O_K}} \arrow[d] \\
    {\RG_{\cris}((R/p) / W)} \arrow[r] & \RG_{\cris}(R_0/W) \arrow[r] & {\RG_{\cris}(R_0/W) \otimes_{W} \mathcal O_K \simeq \RG_{\cris}(R_0/\mathcal O_K)}
\end{tikzcd}\]
Here $\mathcal O_K$ is given the (unique) PD structure on the maximal ideal.

The natural transformation $\chi_i \colon \Q_p(i)(-/p) \to \RG_{\cris}((-/p) / \Z_p)$ of functors on $\mrm{qSyn}_{\Z_p}$ comes (via composition with $(-)/p$) from a natural transformation \[\tilde{\chi}_i \colon \Q_p(i)(-) \to \RG_{\cris}(-/\Q_p)\] of functors on $\mrm{qSyn}_{\F_p}$ by \cite[Corollary 6.20]{AMMN}. Hence the composite \[\Q_p(i)(R/p) \xrightarrow{\chi_i} (\dR_{R/\Z_p})_{\Q_p} \to \RG_{\cris}(R_0/\Q_p)\] identifies (functorially in $R$) with the composite \[\Q_p(i)(R/p) \xrightarrow{\simeq} \Q_p(i)(R_0) \xrightarrow{\chi_i} \RG_{\cris}(R_0 / \Q_p).\]
We now observe that for any $S \in \mrm{qSyn}_{\F_p}$, $\chi_i \colon \Q_p(i)(S) \to \RG_{\cris}(S/\Q_p)$ factors as $\Q_p(i)(S) \xrightarrow{\chi_i'} \RG_{\cris}(S/\Q_p)^{\phi = p^i} \to \RG_{\cris}(S/\Q_p)$. This map $\chi_i'$ is an equivalence by quasisyntomic descent and \cite[Proposition 6.18]{AMMN}.

Putting this together, we find a functorial homotopy as in the statement.
\end{proof}

\begin{proof}[Proof of \cref{thmA}]

We first show that (1) is equivalent to (2). Let us assume temporarily that $k$ is algebraically closed. By \cite[Theorem 7.5.6]{BL}, there is an equivalence $\RG_{\et}(-, \G_m^{\land}) \xrightarrow{\simeq} \RG_{\syn}(-, \Z_p(1))[1]$ of functors on $p$-adic formal schemes. By \cref{simplify}, $\mathfrak X$ is a quasicanonical lift of $\mathfrak X_0$ if and only if the map $\mrm{H}^d_{\mrm{syn}}(\mathfrak X, \Q_p(1)) \to \mrm{H}^d_{\mrm{syn}}(\mathfrak X_0, \Q_p(1))$ is surjective. Considering \cref{var}, this is the case if and only if the image of $\mrm{H}^d_{\mrm{syn}}(\mathfrak X_0, \Q_p(1)) \to \mrm{H}^d_{\dR}(\mathfrak X/\mathcal O_K)[1/p]$ is contained in that of $\mrm{H}^d(\mathfrak X, \Omega^{\ge 1}_{\mathfrak X/\mathcal O_K})[1/p] \to \mrm{H}^d_{\dR}(\mathfrak X/\mathcal O_K)[1/p]$. Using the identification $\mrm{H}^d_{\dR}(\mathfrak X_K/K) \simeq \mrm{H}^d_{\cris}(\mathfrak X_0/K_0) \otimes_{K_0} K$, this is precisely the condition that \[\mrm{H}^d_{\cris}(\mathfrak X_0/K_0)^{\phi = p} \subset \mrm{Fil}^1 \mrm{H}^d_{\dR}(\mathfrak X_K/K).\] Since the slope $1$ piece $\mrm{H}^d_{\cris}(\mathfrak X_0/K_0)^{(1)} \subset \mrm{H}^d_{\cris}(\mathfrak X_0/K_0)$ is spanned as a $K_0$-module by $\mrm{H}^d_{\cris}(\mathfrak X_0/K_0)^{\phi = p}$, the condition is equivalent to $\mrm{H}^d_{\cris}(\mathfrak X_0/K_0)^{(1)} \otimes_{K_0} K \subset \mrm{Fil}^1 \mrm{H}^d_{\dR}(\mathfrak X_K/K)$.

More generally, if $k$ is not necessarily algebraically closed, then we can still check whether $\mrm{H}^d_{\cris}(\mathfrak X_0/K_0)^{(1)} \otimes_{K_0} K$ is contained in $\mrm{Fil}^1 \mrm{H}^d_{\dR}(\mathfrak X_K/K)$ after base change along $K \to \breve{K}$, whereafter it becomes the condition that $\mrm{H}^d_{\cris}(\overline{\mathfrak X}_0/ \breve{K}_0)^{(1)} \otimes_{\breve{K}_0} \breve{K} \subset \mrm{Fil}^1 \mrm{H}^d_{\dR}(\mathfrak X_{\breve{K}}/\breve{K})$. By the above, this corresponds to $\mathfrak X_{\mrm{Spf}(\breve{\mathcal O}_K)}$ being a quasicanonical lift of $\overline{\mathfrak X}_0$, or, in view of \cref{simplify} (2), $\mathfrak X$ being a quasicanonical lift of $\mathfrak X_0$. 

We now show that (2) and (3) are equivalent. Firstly, if (2) holds, then $D = \mrm{H}^d_{\cris}(\mathfrak X_0/K_0)^{(1)}$ equipped with the subspace filtration on $D_K$ defines a weakly admissible sub-(filtered $\phi$-module) of the crystalline cohomology. Indeed, $D$ being a subobject of a weakly admissible sub(-filtered $\phi$-module), it suffices to show that the endpoint $t_H(D) \ge t_N(D)$, which follows as the only Newton slope is $1$ and the Hodge slopes are $\ge 1$. Hence $D$ corresponds to a subrepresentation $V(-1) \subset \mrm{H}^d_{\et}(\mathfrak X_C, \Q_p)$ with $V$ unramified. Moreover, since the quotient of $\mrm{H}^d_{\cris}(\mathfrak X_0/K_0)$ by $D$ has no slope $1$ piece, the corresponding Galois representation has no Hodge--Tate weight $1$ piece by \cref{lemm} below.

Conversely, if (3) holds, then there exists a weakly admissible sub filtered $\phi$-module $D$ of pure Hodge slope $1$ such that the quotient has no Hodge slope $1$. By another application of \cref{lemm}, we see that $D$ is $\mrm{H}^d_{\cris}(\mathfrak X_0/K_0)^{(1)}$ with the subspace filtration, as is desired.
\end{proof}

\begin{lemma} \label{lemm} Let $Y$ be a Bloch--Kato $n$-ordinary smooth and proper $k$-scheme, and suppose $\mathfrak X/\mrm{Spf}(\mathcal O_K)$ is a smooth lift of $Y$ with rigid-analytic generic fibre $V$.

Then the only Newton and Hodge slopes of the filtered $\phi$-modules attached to $\mathfrak X$ in the interval $[0,n)$ are integers, and the first $n$ segments of the Hodge polygon of the generic fibre and Newton polygon of the special fibre agree. That is, for each $i < n$ and every $j$, \begin{equation} \label{equal} \dim_K \mrm{H}^{j}(\mathfrak X_K, \Omega^i_{\mathfrak X_K/K}) = \dim_{K_0} \mrm{H}^{i+j}_{\cris}(Y/K_0)^{(i)}.
\end{equation}
\end{lemma}

\begin{proof}

We may assume that $k$ is algebraically closed. It is standard that the slope $[i,i+1)$ piece of $\mrm{H}^{m}_{\cris}(Y/K_0)_{\Q}$ is given by $\mrm{H}^{m-i}(W\Omega^i_Y)_{\Q}$ \cite[Corollaire 3.5.]{dRW-Illusie}. For $Y$, these are pure of slopes $i$ for each $i < n$ by (the proof of) \cite[Proposition 7.3]{BlochKato}, since 
 \begin{equation} \label{BK2-iso} \mrm{H}^{*}(Y,W\Omega^i_Y) \cong \mrm{H}^*_{\et}(Y, W\Omega_{\mrm{log}}^i) \otimes_{\Z_p} W(k).
 \end{equation}
 It now remains to prove that $\dim_K \mrm{H}^{j}(V, \Omega^i_{V/K}) = \dim_{K_0} \mrm{H}^{j}(Y, W\Omega^i)[1/p]$ for each $i < n$ and every $j$.

We first sketch a proof under additional torsion assumptions; this suffices for \cref{thmA} under the additional hypothesis that $\mrm{H}^{d-1}_{\et}(\mathfrak X_0, \mathbf G_m^{\land})$ is torsionfree. By examining the proof of \cite[Proposition 7.3]{BlochKato}, BK $n$-ordinarity implies that $\RG(Y, W\Omega^i)/p \simeq \RG(Y, \Omega^i)$ for each $i < n$. For each $j$, we now have a string of inequalities
\[  \dim_{K} \mrm{H}^{j-i}(\Omega^i_{V/K}) \le \dim_k \mrm{H}^{j-i}(\Omega^i_{Y/k}) =\dim_k \mrm{H}^{j-i}_{\et}(W\Omega^i_{Y/k})/p  = \dim_{K_0} \mrm{H}^{j-i}_{\et}(W\Omega^i)[1/p],\]
provided we assume that $\mrm{H}^{j-i}(W\Omega^i_{Y/k}), \mrm{H}^{j-i+1}(W\Omega^i_{Y/k})$ are torsion free; cf.\ the proof of \cite[Lemma 7.1]{BlochKato}. Since the Newton polygon lies above the Hodge polygon, we also have the reverse inequalities, so we may conclude.

More generally, we can instead appeal to a recent result of Y.\ Yang relating the Hodge numbers $h^{a,b}$ of the generic fibre $V$ to the Hodge--Witt numbers $h_W^{a,b}$ of the special fibre; we refer to \cite[\S 2]{yyang26} for the notation used in the rest of this proof. The $W$-modules $\mrm{H}^*(Y, W\Omega^a)$ is finitely generated for each $a < n$ (e.g.\ by the proof of \cite[Proposition 7.3]{BlochKato}), so by \cite[II, Corollaire 3.9]{IllusieRaynaud}, the domino numbers $T^{a,*}$ vanish in this regime. Hence for each $a < n$, $h^{a, b-a}_W$ coincides with the slope number $m^{a, b-a}$, which, by the earlier statement about slopes (and the definition \cite[(2.1)]{yyang26}) is simply $\dim_{K_0} \mrm{H}^{b}_{\cris}(Y/K_0)^{(a)}$. Now applying \cite[Theorem 1.1]{yyang26}, for each $r \le n$ and every $b$ we have
\[\sum_{a < r} (r-a) h^{a,b-a} \le \sum_{a < r} (r-a) \dim_{K_0} \mrm{H}^{b}_{\cris}(Y/K_0)^{(a)},\]
i.e.\ the first $n$ segments of the generic Hodge polygon lie on or above those of the Newton polygon of the special fibre. Since the converse always holds (for smooth proper formal schemes), all of these inequalities are equalities, and the result follows.
\end{proof}

\section{Complements}

We now establish some additional properties of quasicanonical lifts and compatibilities.

\subsection{Algebraisability}
\begin{proposition} Let $Y$ be a BK $2$-ordinary Calabi--Yau variety over $k$ such that $\mrm{H}^d_{\et}(Y_{\overline{k}}, \Z_p)$ is torsionfree.

 \begin{enumerate} \item  If $\mathscr L$ is a line bundle over $Y$ and $X$ is a flat lift of $Y$ to $A \in \mrm{CRing}^{\mrm{art}}_W$, corresponding to an $m$-torsion class $c \in \pi_0 \mrm{ST}_Y(A)$, then $\mathscr L^{\otimes m}$ lifts to $X$. 
 \item Suppose moreover that $Y$ is projective and $\mrm{H}^1(Y, \mathcal O_Y) = 0$. Then for any complete local Noetherian ring $R$ with residue field $k$, every formal quasicanonical lift $\mathfrak X/\mrm{Spf}(R)$ algebraises to a projective $R$-scheme.
 \end{enumerate}
\end{proposition}

The analogue for K3 surfaces was proven by Nygaard and Yu \cite[Proposition 1.8]{NygaardTC} \cite[Lemma 1.4]{JDY}.

\begin{proof}

We first show how (2) follows from (1). The hypothesis that $h^1(\mathcal O_Y)=0$ guarantees that lifts of line bundles are unique up to isomorphism if they exist, so the lifts of $\mathscr L^{\otimes m}$ to each Artinian ring in (1) are necessarily mutually compatible. Hence $\mathscr L^{\otimes m}$ lifts to a line bundle on $\mathfrak X$. Taking $\mathscr L$ to be ample, this shows that $\mathfrak X$ is algebraisable by the Grothendieck existence theorem \cite[Theorem 089A]{stacks-project}.
   
We now prove (1). It is convenient to place this inside Brantner--Taelman's formalism \cite[\S 8]{BT}, which we briefly recall for the reader's convenience. Recall that a line bundle over a scheme $Z$ corresponds to a map $\lambda \colon \Z[-1] \to \G_m$ of sheaves on $Z$. Given a line bundle $\mathscr L$ over $Y$ (corresponding to a map $\lambda$ as above), Brantner--Taelman construct a deformation functor $\mrm{Def}_{Y, \lambda}$ whose $A$-points parametrise lifts of the pair $(Y, \mathcal L)$, and a variant $\mrm{ST}_{Y, \lambda}$ of the Serre--Tate period domain fitting into a commutative diagram of formal moduli problems
  \[\begin{tikzcd} {\mrm{Def}_{Y, \lambda}(-)} \arrow[r, "\simeq"] \arrow[d, "\text{forget}"] & {\mrm{ST}_{Y, \lambda} = \Map_{\mrm{Mod}(S_{\et}, \Z_p)} (\mrm{R}f_* P_{Y, \lambda}, \mrm{R}^d f_* \mathbf Z_p \otimes_{\Z_p} \mathbf G_m^{\land}(\mathcal W \otimes_W (-))[1-d])} \arrow[d, "\varphi"] \\
      {\mrm{Def}_{Y}(-)} \arrow[r, "\simeq"] & {\mrm{ST}_{Y} = \Map_{\mrm{Mod}(S_{\et}, \Z_p)} (\mrm{R}f_* \G_m^{\land}, \mrm{R}^d f_* \mathbf Z_p \otimes_{\Z_p} \mathbf G_m^{\land}(\mathcal W \otimes_W (-))[1-d]).}
  \end{tikzcd}\]
Here $P_{Y,\lambda} = \mrm{cofib}(\lambda \colon \Z_p[-1] \to \mathbf G_{m,Y}^{\land})$, $\varphi$ is given by composition with $\mrm{R}f_* \mathbf G_m^{\land} \to \mrm{R}f_* P_{Y,\lambda}$, and (as usual) we have written $f$ for the structure morphism $Y \to \Spec k$.

Now suppose we are in the situation of (1). Take $\mathscr L$ on $Y$, corresponding to $\lambda \colon \Z[-1] \to \G_m$, so $\mathscr L^{\otimes m}$ corresponds to $\mu = m \lambda$ (writing $\mathbf G_m$ additively). Then we can consider the exact sequence
\[\pi_0 \mrm{ST}_{Y, \mu}(A) \xrightarrow{\varphi} \pi_0 \mrm{ST}_Y(A) \xrightarrow{\phi'} \pi_0 \mrm{Map}_{\mrm{Mod}(S_{\et}, \Z_p)}(\mrm{R}f_* \Z_p[-1], \mrm{R}^d f_* \mathbf Z_p \otimes_{\Z_p} \mathbf G_m^{\land}(\mathcal W \otimes_W A)[1-d]).\]
The image of $c$ under $\phi'$ vanishes: it is the homotopy class of the composite $\mrm{R}f_* \Z_p[-1] \xrightarrow{m \lambda} \mrm{R}f_* \mathbf G_m^{\land} \xrightarrow{c} \mrm{R}^d f_* \mathbf \Z_p \otimes_{\Z_p} \G_m^{\land}(\mathcal W \otimes_W A)[1-d]$, and $mc = 0$ by hypothesis. It follows that $c$ lies in the image of $\pi_0 \mrm{ST}_{Y, \mu}(A)$, i.e.\ that $\mathscr L^{\otimes m}$ lifts to $X$.
\end{proof}

\subsection{Canonical lifts and compatibilities}

We now consider some compatibilities with other (quasi)canonical lifts in the literature using our characterisation of quasicanonical lifts, and show that in some cases it is possible to single out the canonical lift among the quasicanonical lifts. 

By \cite[Corollary 5.1.5, Remark 5.1.6]{AZ-ST}, every $1$-ordinary Calabi--Yau variety over $k$ admits a canonical lift to $W_2(k)$, constructed explicitly in terms of the (unique) $F$-splitting.\footnote{Recall that an $F$-splitting of $Y$ is an $\mathcal O_Y$-linear retraction of the Frobenius $F \colon \mathcal O_Y \to F_* \mathcal O_Y$.} Achinger--Zdanowicz develop a variant of Serre--Tate theory in \cite[\S 11]{AZ-ST} under more restrictive hypotheses, and in particular produce a canonical lift mod $p^2$ of any BK $2$-ordinary Calabi--Yau variety, which they show to agree with the canonical lift coming from the $F$-splitting \cite[Proposition 11.3.1]{AZ-ST}. Their argument readily generalises to give point (1) below.
\begin{proposition} \label{BTAZ} Let $Y$ be a BK $2$-ordinary Calabi--Yau variety over $k$ such that $\mrm{H}^d_{\et}(Y_{\overline{k}}, \Z_p)$ is torsionfree. 
  \begin{enumerate} \item Brantner--Taelman's canonical lift to $W_2(k)$ is isomorphic to that of Achinger--Zdanowicz via $F$-splittings.
    \item If $Y$ satisfies $(\dagger)$, then the canonical formal lift $\mathfrak X/\mrm{Spf}(W)$ is the unique \emph{quasicanonical} lift which agrees mod $p^2$ with Achinger--Zdanowicz' lift. 
    \end{enumerate}
  \end{proposition}

\begin{remark} Of course, part (2) is only relevant when $p =2$, so that $W(\overline{k})$ possesses a nontrivial $p$th root of unity.
\end{remark}

  \begin{proof}[Proof of \cref{BTAZ}] For (1), let $Z$ denote Achinger--Zdanowicz's lift to $W_2(k)$. The map $\mathbf G_m(\mathcal O_{Z}) \to \mathbf G_m(\mathcal O_Y)$ admits a multiplicative section \cite[Lemma 5.2.2]{AZ-ST}; hence the map $\RG(Z_{W_2(\overline{k})}, \mathbf G_m^{\land}) \to \RG(Y_{\overline{k}}, \G_m^{\land})$ admits a section and the map to the cofibre vanishes; this means $Z$ is Brantner--Taelman's lift. For part (2), it is enough, in view of (1), to show that a quasicanonical lift to $\mrm{Spf}(W(k))$ is the canonical lift if and only if it is so mod $p^2$. There is nothing to prove unless $p=2$, in which case, by \cref{simplify} (1), this reduces to the observation that $1$ and $-1$ are distinct mod $4$.
\end{proof}

\cref{BTAZ} allows us to check some compatibilities:

\begin{proposition} Let $Y/k$ be as in \cref{sit}. If $Y$ has a lift $\mathfrak X/\mrm{Spf}(W)$ to which the Frobenius lifts, then $\mathfrak X$ is the canonical lift of $Y$. In particular, if $Y$ is an ordinary abelian variety, then $Y^{\mrm{can}}$ agrees with the lift due to Serre--Tate.
\end{proposition}

\begin{proof} 
We recall that if $(\mathfrak X,F)$ is a lift of $Y$ with its Frobenius, then, under the equivalence $\RG_{\dR}(\mathfrak X/W) \simeq \RG_{\cris}(Y/W)$, the crystalline Frobenius $\phi$ corresponds to that induced on de Rham cohomology by $F$; it follows that the Hodge filtration is $\phi$-stable. Hence $D = \mrm{Fil}^1 \mrm{H}^d_{\cris}(\mathfrak X/K_0)$ with the restricted Frobenius defines a sub-(filtered $\phi$-module), and by considering its Newton polygon, we must have $\mrm{H}^d_{\cris}(Y/K_0)^{(1)} \subset D$. Hence \cref{thmA}, $\mathfrak X$ is a quasicanonical lift of $Y$. Moreover, Achinger--Witaszek--Zdanowicz prove the analogue mod $p^2$ \cite[Theorem 2.7]{AWZ2}. Therefore \cref{BTAZ} (2) implies that $\mathfrak X$ is the canonical lift of $Y$. If $Y$ is an ordinary abelian variety, then condition $(\dagger)$ automatically holds -- indeed, all the crystalline cohomology is torsionfree -- and the Serre--Tate canonical lift to $W(k)$ is uniquely characterised by its lifting of Frobenius \cite[Appendix, Corollary 1.2]{Messing}.
\end{proof}

\section{Application to derived invariance} \label{addendum}

We now present a quick application of our characterisation to derived categories using the author's \cite{moore1} and results towards derived invariance of rational motives. This is an at least preliminary higher-dimensional analogue of \cite[Theorem 1.1]{kuramaFM} (which concerns ordinary abelian varieties and K3 surface) and \cite[Proposition 4.5.15]{Ward} (concerning threefolds at $p > 3$). It also partly answers \cite[Question 1.9]{kuramaFM}. For our purposes, we take a \emph{derived equivalence} between $A$-schemes to be an equivalence of $A$-linear $\infty$-categories of quasicoherent sheaves; these are equivalently (under the present conditions) Fourier--Mukai equivalences.

\begin{proposition} \label{der-equiv} Suppose $Y, Y'$ are two projective BK $2$-ordinary Calabi--Yau varieties over a perfect field $k$ of characteristic $p$ satisfying the following conditions (where $Z = Y,Y'$):
\begin{enumerate}
 \item The conditions of \cref{thmA}; i.e.\ $\mrm{H}^d_{\et}(\overline{Z}, \Z_p)$ and $\mrm{H}^d_{\et}(\overline{Z}, \G_m^{\land})$ are torsion free
\item The groups $\mrm{H}^1(Z, \mathcal O_Z), \mrm{H}^2(Z, \mathcal O_Z)$, and $\mrm{H}^0(Z, \Omega^{d-2})$ vanish.
\end{enumerate}

If $p > 2$, then any derived equivalence between $Y$ and $Y'$ lifts to a derived equivalence between the canonical lifts $X, X'$.

More generally, for any derived equivalence between $Y$ and $Y'$ and quasicanonical lift $X$ of $Y$ to $W(k)$, there exists a unique quasicanonical lift $X'$ of $Y'$ such that there exists a derived equivalence between $X$ and $X'$ lifting the specified derived equivalence.

\end{proposition}

\begin{remark}
  \begin{enumerate} \item For example, we can apply \cref{der-equiv} to strict Calabi--Yau variety (necessarily of dimension $\ge 3$) with mild Hodge symmetry (in turn implied by Hodge--de Rham spectral sequence degeneration and having torsionfree crystalline cohomology).
  \item The derived equivalence generally does not lift uniquely as the objects can have nontrivial automorphisms over $W$. In particular, this is why the phrasing at $p=2$ is more delicate.
  \end{enumerate}
\end{remark}

\begin{proof}[Proof of \cref{der-equiv}]
  
Let $X$ be a quasicanonical lift of $Y$ to $W = W(k)$. By \cite[Corollary 1.4]{moore1}, taking the derived category induces an isomorphism between the (classical) deformation functors of each $Z = Y,Y'$ and their derived categories.\footnote{Note that there is no issue with `completing' the formal moduli problem of categorical deformations in this case.} Hence there exists a unique formal lift $\mathfrak X'/\mrm{Spf}(W)$ of $Y'$ and a coherent collection of equivalences $\QC(\mathfrak X'/p^n) \simeq \QC(X/p^n)$. Since $\mrm{H}^2(Y', \mathcal O)$, any ample line bundle on $Y'$ lifts to $\mathfrak X'$; by the Grothendieck algebraisation theorem \cite[Theorem 089A]{stacks-project}, $\mathfrak X'$ algebraises (uniquely) to a projective $W$-scheme $X'$. By the derived Grothendieck existence theorem \cite[Theorem 0DIG]{stacks-project} applied to the (perfect) corresponding Fourier--Mukai kernel on $\mathfrak X \times_{\mrm{Spf}(W)} \mathfrak X'$, the derived equivalence lifts to one between $X$ and $X'$. Passing to generic fibres, we obtain a derived equivalence between the generic fibres $V ,V'$. It remains to show that $X'$ is a quasicanonical lift of $Y'$.

  Write $d = 2a+\epsilon$ where $\epsilon \in \{0,1\}$ and put $U = \bigoplus_{m} \mrm{H}^{2m+\epsilon}_{\et}(\overline{V}, \Q_p)(m-a)$ (and analogously for $U')$. By \cite[Lemma 3.1]{Honigs} or \cite[\S 2]{LieblichOlsson}, the derived equivalence between $V,V'$ induces an isomorphism $U \cong U'$ of $G_{K}$-representations. Let $D_{m}$ be the filtered $\phi$-module $\mrm{D}_{\cris}\Big(\mrm{H}^{2m+\epsilon}_{\et}(\overline{V}, Q_p)\Big)$, i.e. $\mrm{H}^{2m+\epsilon}_{\cris}(Y/K)$ equipped with the Hodge filtration from $\mrm{H}^{2m+\epsilon}_{\dR}(V/K)$. $D = D_{\cris}(U)$ is given by $\bigoplus_m D_{m}(m-a)$.

  If $m \ne a$, then we automatically have $D_m(m-a)^{(1)} \subset \mrm{Fil}^1(D_m(m-a))$, as in fact $\mrm{Fil}^{1+m-a} \mrm{H}^{2m+\epsilon}_{\dR}(V/K) = \mrm{H}^{2m+\epsilon}_{\dR}(V/K)$. If $m < a$ this is trivial; if $m > a$, we observe that if $i + j = 2m+ \epsilon$ and the Hodge number $h^{i,j}(V) \ne 0$, then necessarily $i \ge 2m + \epsilon - \dim(Y) = 2(m-a) \ge 1 + m-a$. Hence the condition that $\mrm{H}^d_{\cris}(Y/K)^{(1)} \subset \Fil^1 \mrm{H}^d_{\dR}(V/K)$ is equivalent to $D^{(1)} \subset \Fil^1 D$. The same conclusion holds for $Y'$ and the (analogously defined) filtered $\phi$-module $D'$. By \cref{thmA}, and since $D \cong D'$, we conclude that $X'$ is a (quasi)canonical lift of $Y'$.
\end{proof}

\section{A canonical lift not of CM-type} \label{noncm}

We now move on to study a canonical lift of an ordinary Calabi--Yau variety which possesses some peculiar properties. The example comes from the (strict) projective Calabi--Yau variety $V$ over $\Q$ appearing in \cite{Saito-Yui}, where it is called $\tilde{S}^2_{\Gamma_1(6)}$. This has a smooth projective model $X \to \Spec \Z[1/6]$, and its mod $p$ reductions will produce the examples in \cref{thmB}.

\subsection{Recollections on $X$. }

We first recall some properties of the variety $X/\Z[1/6]$ following \cite[\S 4]{Saito-Yui} and \cite[\S 1 and Lemma 3.1]{Schoen88}.

The details of its construction are mostly unimportant for our purposes, as we will focus on its cohomology, but let us give a brief overview. Consider the semistable elliptic surface $S_1(6) \to \mathbf P^1_{\Z}$ \cite[Theor\`eme]{Beauville}; blowing up the three nodes produces a smooth projective surface $S_1(6)$ over $\Z$. Restricting to $R = \Z[1/6]$, the $R$-scheme $Z = S_1(6)_{R} \times_{\mathbf P^1_R} S_1(6)_{R}$ is smooth away from the (geometric) points lying in a common singular fibre over $\mathbf P^1$, which constitute $50$ ordinary double points. By blowing up along smooth surfaces containing each node in some order, we obtain a small projective resolution $X \to Z$ \cite[Lemma 3.1]{Schoen88}. 

\begin{definition} We refer to $X \to \Spec \Z[1/6]$ as the integral model; we write $V = X_{\Q}$ for the generic fibre, and put $X_p = X_{\Spec{\Z_p}}$ for $p \ge 5$.
 \end{definition}

By \cite[Corollary 4.2]{Saito-Yui}, $V$ is a rigid strict Calabi--Yau threefold. 

 \subsubsection{Galois representation}
 
The following modularity property was observed by Saito--Yui, who attribute it to Kuga--Sato--Shimura and Deligne \cite[Theorem 4.1]{Saito-Yui}.  We remark that now \emph{all} rigid strict Calabi--Yau threefolds over $\Q$ are known to be modular by a theorem of Gouv\^{e}a--Yui \cite[Theorem 3]{GouvYui}, and one can even identify $f$ from their result in this case.

\begin{theorem}\label{modular} Let $f$ be the newform $f(z) = \eta(z)^2 \eta(2z)^2 \eta(3z)^2 \eta(6z)^2$ of weight $4$ and level $6$. For any prime $\ell$, there is an isomorphism $\mathrm{H}^3_{\et}(V_{\overline{\Q}}, \Q_\ell)^{\vee} \simeq \rho_{f, \ell}$ of $G_{\Q}$-representations.
\end{theorem}

\subsection{The $\ell$-adic monodromy group} \label{cm}

\begin{definition} Let $F$ be a finitely generated field extension of $\mathbf Q$ and $Z$ be a smooth projective variety over $F$ of pure dimension $d$. Fix an algebraic closure $\overline{F}/F$ and a prime $\ell$. Then the \emph{$\ell$-adic algebraic monodromy group} $G_\ell(Z)$ is the algebraic group defined as the identity component of the $\Q_\ell$-Zariski closure of the image of the $\ell$-adic representation $G_{F} \to \mrm{GL}(\mrm{H}^d_{\et}(Z_{\overline{F}}, \Q_\ell)).$
\end{definition}

The group $G_\ell$ is insensitive to shrinking the field of definition \cite[Remark 2.2.2 (i)]{Moonen17}. Hence, given a variety over $\C$, we can (canonically) define its $\ell$-adic algebraic monodromy group as that of any model over a subfield of $\C$ which is finitely generated over $\Q$.

The $\ell$-adic monodromy group of $V$ can be computed using the following result of Ribet. We state it only for even $k$, as $\mathfrak a_\ell$ is defined differently when $k$ is odd.

\begin{proposition}\label{ribet} \emph{\cite[Theorem 5.7]{ribetgalois}} Let $f$ be a non-CM newform of even weight $k \ge 2$ and level $N$, $K = \Q( a_p(f) : p \nmid N)$ and $L \subset K$ be the maximal totally real subfield. Let $G_\ell(f) \subset \mrm{GL}_2(K \otimes_\Q \Q_\ell)$ be the image of the corresponding $\ell$-adic representation $\rho_{f,\ell}$, let $\mathfrak g_\ell$ be its Lie algebra, and put $\mathfrak a = \{u \in \mathfrak{gl}_2(L) \mid \mrm{tr}(u) \in \Q\}$. Then $\mathfrak g_\ell(f) \subseteq \mathfrak a \otimes_\Q \mathbf Q_\ell$, where the inclusion is strict if and only if $f$ admits an extra twist in the sense of \emph{[Ibid, \S5]}.
\end{proposition}

\begin{proposition} \label{monodromy} The $\ell$-adic monodromy group $G_\ell(V)$ is $\mrm{GL}(\mrm{H}^3_{\et}(V_{\overline{\Q}}, \Q_\ell))$. 
\end{proposition}
\begin{proof}
It suffices to prove this on the level of Lie algebras. The newform $f$ attached to $V$ has $k=4$ and $K = L = \Q$; moreover, $f$ does not have CM, and it cannot admit extra twists since $K = \Q$ \cite[modular form 6.4.a.a]{lmfdb}. The result therefore follows immediately from \cref{ribet}.
\end{proof}

\begin{remark} 
  According to the Mumford--Tate conjecture (MTC), for any smooth and proper complex algebraic variety $V$ and prime $\ell$, there should exist an isomorphism $G_\ell(V) \simeq \mrm{MT}(V) \otimes_{\Q} \Q_\ell$, where $\mrm{MT}(V)$ denotes the Mumford--Tate group. The notion of CM type is often formulated in terms of $\mrm{MT}(V)$ being a torus, which is equivalent to the above-established abelianity provided that $V$ satisfies the MTC. We note that, whilst this conjecture remains open in general, it has been proven for K3 surfaces \cite{Tankeev} and abelian varieties of CM type \cite[Theorem 5]{pohlmann}.
\end{remark}
\subsection{$p$-adic invariants}

We now shift our attention to the reductions $Y_p = \mathfrak X_{\F_p}$ for $p \ge 5$.
\begin{proposition} \label{cy-check} For each $p \ge 5$, $Y_p$ is a Calabi--Yau variety over $\F_p$.
\end{proposition}

\begin{proof} We need to show that $Y_p$ has trivial canonical bundle and is geometrically irreducible. Since $\mathfrak X_p$ has geometrically irreducible generic fibre $V_{\Q_p}$, its special fibre $Y_p$ is geometrically irreducible by \cite[Lemma 0AY8]{stacks-project}; this can also be deduced from the calculation of $h^0(Y_p, \mathcal O) $ below. Moreover, since $V$ has trivial canonical bundle, $\Omega^3_{\mathfrak X_p/\Z_p}$ restricts to the trivial bundle over the generic fibre, and is therefore itself trivial; this is hence true of the special fibre too by base change. 
\end{proof}
\begin{proposition} \label{modp}\, Let $p \ge 5$ be a prime. Then 
  \begin{enumerate} \item The groups $\mrm{H}^*_{\et}(V_{\overline{\Q}}, \Z_p)$ are torsionfree.
  \item The crystalline cohomology $\mrm{H}^*_{\cris}(Y_p/\mathbf Z_p)$ is torsionfree.
  \item The Hodge numbers of $Y_p$ coincide with those of $V$. In particular, $Y_p$ is a rigid strict Calabi--Yau variety.
  \end{enumerate}
\end{proposition}

\begin{proof}
  The first part can be read off from \cite{Schoen13} with some care. Recall that $t \colon V_{\overline{\Q}} \to Z_{\overline{\Q}}$ is a small resolution (which we will write as $\overline{V} \to \overline{Z}$ for convenience), whose exceptional fibres are projective lines. Schoen considers a `big' resolution $W \to \overline{Z}$ obtained by blowing up the reduced singular locus of $\overline{Z}$, and \cite[Theorem 0.3]{Schoen13} implies that $\mrm{H}^*_{\et}(W, \Z_\ell)$ is torsionfree. But by (e.g.) the description of \cite[p.\ 1, end of para.\ 1]{Schoen88}, $W$ is a blow up of $\overline{V}$ along the exceptional $\mathbf P^1$s, and so the blowup formula implies that $\mrm{H}^*_{\et}(\overline{V}, \Z_\ell)$ is also torsionfree.\footnote{Alternatively, Schoen proves that $\mrm{H}^i_{\et}(\overline{Z}, \Z_\ell)$ is torsionfree for $i \ne 4$ \cite[Proposition 4.1]{Schoen13}; one can use the map $\overline{V} \to \overline{Z}$ and Poincar\'{e} duality to provide an alternative proof of (1).}
    A result of Breuil \cite[Corollaire 4.3.1.5]{Breuil99} implies that for each $i \le p-2$, the $\Z_p$-modules $\mrm{H}^i_{\et}(V_{\overline{\Q}_p},\Z_p)$ and $\mrm{H}^i_{\dR}(\mathfrak X_p/\Z_p) \simeq \mrm{H}^i_{\cris}(Y_p/\Z_p)$ have the same invariant factors. This shows that (1) implies (2). To accomodate the cases $p = 5,7$, we can use use the torsion pairing on integral crystalline cohomology.

Since $p \ge 3 = \dim(Y_p)$ and $Y_p$ lifts to $\Z/p^2$, the Hodge--de Rham spectral sequence degenerates at $\mrm{E}_1$ to give $\dim_{\F_p} \mrm{H}^*_{\dR}(Y_p/\F_p) = \sum_{i,j} h^{i,j}_{Y_p}$, and similarly $\dim_{\Q_p} \mrm{H}^*_{\dR}(V_{\Q_p}/\Q_p) = \sum_{i,j} h^{i,j}_{V}$ \cite[Corollaire 2.4]{DI}. As the crystalline cohomology is torsion free, the dimensions $\dim_{\F_p} \mrm{H}^*_{\dR}(Y_p/\F_p)$ and  $\dim_{\Q_p} \mrm{H}^*_{\dR}(V_{\Q_p}/\Q_p)= \dim_{\mathbf Q} \mrm{H}^*_{\dR}(V/\Q)$ coincide; since we always have $h^{i,j}_{V} \le h^{i,j}_{Y_p}$ for each $(i,j)$, we must have equality throughout. Combined with \cref{cy-check}, we can now conclude that $Y_p$ is a rigid strict Calabi--Yau variety.
\end{proof}

\subsection{Ordinarity}

We finally show that the mod $p$ reduction $Y_p$ is ordinary for many choices of $p$. More precisely, we claim: 
\begin{proposition} \label{iso} For each $p \ge 5$ and $i \ne 3$, the isocrystal $\mrm{H}^i_{\cris}(Y_p/\Q_p)$ is ordinary (that is, its Newton and Hodge polygons coincide). If $5 \le p < 47351$, then so is $\mrm{H}^3_{\cris}(Y_p/\Q_p)$.
\end{proposition}

\begin{remark} This definition of ordinarity (in terms of the isocrystals) is equivalent to various others thanks to the absence of crystalline torsion; see \cite[Proposition 7.3]{BlochKato}. 
\end{remark}
\begin{proof}[Proof of \cref{iso}]
  
By Poincar\'{e} duality, we may assume $0 \le i \le 3$. The case $i = 0$ is trivial, and $\mrm{H}^1_{\cris}(Y_p/\Q_p)=0$ by \cref{modp}. For $i=2$, the filtered $\phi$-module attached to $\mathfrak X_p$ is pure of Hodge slope $1$, so, since it is weakly admissible, ordinarity is immediate.

Now consider $i = 3$, and fix a prime $\ell \nmid 6p$. By \cite[Theorem 1]{Riemann}, the characteristic polynomial of the Frobenius on $\mrm{H}^3_{\cris}(Y_p/\Q_p)$ coincides with that of the geometric Frobenius on $\mrm{H}^3_{\et}((Y_p)_{\overline{\F}_p}, \Q_\ell) \simeq \mrm{H}^3_{\et}(V_{\overline{\Q}_p}, \Q_\ell) = \rho_{f,\ell}^{\vee}|_{\mrm{Gal}(\Q_p)}$, namely $X^2 - a_p(f) X + p^3$. If $p \nmid a_p(f)$, the roots of this polynomial have $p$-adic valuations $0$ and $3$; hence the Newton slopes agree with the Hodge slopes and $\mrm{H}^3_{\cris}(Y_p/\Q_p)$ is ordinary. But $p \nmid a_p$ provided $5 \le p < 47351$ by direct computation with \cite{sagemath}.
\end{proof}

\begin{proposition} \label{decomp}
The representation $\mrm{H}^3_{\et}(V_{\overline{\Q}}, \Q_7)^{\vee}|_{G_{\Q_7}} \simeq \rho_{f,7}|_{G_{\Q_7}}$ is indecomposable.
\end{proposition}

\begin{proof}
  We apply a general method in \cite{Ghate}. According to \cite[Proposition 6]{Ghate}, to show that $\rho_{f,7}|_{G_{\Q_7}}$ does not split, it suffices to show that (a) $f$ is $7$-distinguished (in the sense that the semisimplification of $\bar{\rho}_{f,7}|_{G_7}$ splits as a sum of distinct characters), (b) $\bar{\rho}_{f,7}$ is (globally) absolutely irreducible, and (c) $\bar{\rho}_{f,7}|_{G_{\Q_7}}$ does not split.

It is clear that (a) holds: the semisimplification of $\bar{\rho}_{f,7}|_{G_7}$ has the form $\chi \oplus \chi^{-1} \omega_{\mrm{cyc}}^3$, where $\chi$ is an unramified character. For step (b), we claim that the global mod $7$ representation $\bar{\rho}_{f,7}$ is absolutely irreducible. Otherwise, its semisimplification is a sum $\chi \oplus \chi^{-1} \omega_{\mrm{cyc}}^3$ for some character $\chi$, where $\omega_{\mrm{cyc}}$ is the mod $7$ cyclotomic character. If $\chi$ is ramified at some prime $\ell \ne 7$, then (since $\omega_{\mrm{cyc}}$ is only ramified at $7$), $\ell^2$ divides the level of $f$; this contradicts \cite[Proposition 0.1]{livne}, as the prime to $7$ conductor divides the level of $f$, namely $6$. Hence $\chi = \omega_{\mrm{cyc}}^i$ for some $i$, so $\bar{\rho}^{\mrm{ss}}_f \simeq \omega_{\mrm{cyc}}^i \oplus \omega_{\mrm{cyc}}^{3-i}$. Evaluating at (e.g.) $\mrm{Frob}_5$ we find $6 \equiv a_5 \equiv 5^i + 5^{3-i}$ mod $7$, which is impossible. 

We now need to show (c). Suppose for contradiction that $\bar{\rho}_{f,7}|_{G_{\Q_7}}$ splits as a sum of two characters. Then $\bar{\rho}_f$ is tamely ramified at $7$, and hence has a companion form $g$ by a result of Gross and Coleman--Voloch \cite[Theorem 13.10]{Gross} \cite[Theorem 0.1]{ColemanVoloch}. Any companion form to $f$ mod $7$ would have to be of weight $(7+1) - 4 = 4$, the same as $f$, and since $\dim_{\F_7} S_4(\Gamma_1(6), \F_7) = 1$, $f$ mod $7$ would have to be a companion of itself. This is not the case; for example, $a_5 \not \equiv 5^3 a_5 \mod 7$.
\end{proof}
\begin{proof}[Proof of \cref{thmB}]

The claimed properties of $V$ in (1) and (2) follow immediately from \cref{monodromy}, \cref{modp}, and \cref{iso}. Since $Y_p$ is rigid, the generic fibre of any flat lift to a mixed characteristic complete Noetherian local domain is a base change of $X_p[1/p] \simeq V_{\Q_p}$, so its $\ell$-adic monodromy group agrees with that of $V$. Finally, (3) is simply \cref{decomp}.
\end{proof}

\begin{remark} \label{strength} This example $V$ also shows that some reasonable strengthenings of the conditions in \cref{thmA} (2) do not hold. For example, writing $K = W(k)[1/p]$, it is \emph{not} the case that $\mathrm{Fil}^1 \mrm{H}^d_{\dR}(\mathfrak X_{K}/K) \subset \mrm{H}^d_{\dR}(X_{K}/K)$ is $\phi$-stable for the canonical lift $X/W(k)$ of a BK $2$-ordinary Calabi--Yau variety $Y/k$ as in \cref{thmA}, or that $\mrm{H}^d_{\cris}(X/K)^{(i)} \subset \mrm{Fil}^i \mrm{H}^d_{\dR}(X_K/K)$ for all $i$; this answers \cite[Question 11.3.2]{AZ-ST} in the negative.

  Indeed, were either of these the case for the Schoen threefold $Y_7/\F_7$, then $\mrm{Fil}^1 \mrm{H}^3_{\dR}(V_{Q_7}/\Q_7) = \mrm{Fil}^3 \mrm{H}^3_{\dR}(V_{\Q_7}/\Q_7)$ equipped with the subspace filtration and restricted crystalline Frobenius would define a sub-filtered $\phi$-module $H \subset D = \mrm{H}^3_{\dR}(V_{\Q_7}/\Q_7)$. Such $H$ would be of dimension $1$, pure of Hodge slope $3$, and since $D$ is a weakly admissible ordinary filtered $\phi$-module with slopes $0$ and $3$, $H$ would be of pure Newton slope $3$. Along with the slope $0$ piece, this would split $D$ and hence the corresponding $G_{\Q_7}$-representation $\mrm{H}^3_{\et}(V_{\overline{\Q}_7}, \Q_7)$, contradicting \cref{decomp}.
\end{remark}

\emergencystretch = 1em
\renewcommand*{\bibfont}{\Small}
\printbibliography

\end{document}